
\tolerance=10000
\raggedbottom

\baselineskip=15pt
\parskip=1\jot

\def\sk{\vskip 3\jot}

\def\heading#1{\vskip3\jot{\noindent\bf #1}}
\def\label#1{{\noindent\it #1}}


\def\ref#1;#2;#3;#4;#5.{\item{[#1]} #2,#3,{\it #4},#5.}
\def\refinbook#1;#2;#3;#4;#5;#6.{\item{[#1]} #2, #3, #4, {\it #5},#6.} 
\def\refbook#1;#2;#3;#4.{\item{[#1]} #2,{\it #3},#4.}


\def\({\bigl(}
\def\){\bigr)}

\def\de{\delta}

\def\pq{\{p,q\}}
\def\ineq#1{1/p + 1/q #1 1/2}

{
\pageno=0
\nopagenumbers
\rightline{\tt coxeter.tessellations.tex}
\vskip1in

\centerline{\bf A Census of Vertices by Generations in Regular Tessellations of the Plane}
\vskip0.5in

\centerline{Alice Paul}
\centerline{\tt apaul@hmc.edu}
\sk

\centerline{Nicholas Pippenger}
\centerline{\tt njp@math.hmc.edu}
\sk

\centerline{Department of Mathematics}
\centerline{Harvey Mudd College}
\centerline{1250 Dartmouth Avenue}
\centerline{Claremont, CA 91711}
\vskip0.5in

\noindent{\bf Abstract:}
We consider regular tessellations of the plane as infinite graphs in which $q$ edges and $q$  faces meet at each vertex, and in which $p$ edges and $p$ vertices surround each face.
For $\ineq=$, these are tilings of the Euclidean plane; for $\ineq<$, they are tilings of the hyperbolic plane.
We choose a vertex as the origin, and classify vertices into generations according to their distance (as measured by the number of edges in a shortest path) from the origin.
For all $p\ge 3$ and $q \ge 3$ with $\ineq\le$, we determine the rational generating function giving the number of vertices in each generation.
\vfill\eject
}

\heading{1. Introduction}

A {\it regular tessellation\/} is a planar graph in which every vertex has degree $q\ge 3$ and every face has degree $p\ge 3$.
Following Coxeter [C1], we denote such a graph by $\pq$.
(This notation will {\it not\/} be used to denote a set with two elements.)
When $\ineq>$, the graph $\{p,q\}$ can be drawn on a sphere in a regular way (that is, so that all edges have the same spherical length and all faces the same spherical area).
These tessellations correspond to the Platonic solids:
$\{3,3\}$ is the tetrahedron, $\{4,3\}$ is the cube, $\{3,4\}$ is the octahedron,
$\{5,3\}$ is the dodecahedron, and $\{3,5\}$ is the icosahedron.
When $\ineq=$, the graph $\pq$ can be drawn in the Euclidean plane in a regular way.
These tessellations correspond to tilings of the Euclidean plane by regular polygons:
$\{4,4\}$, $\{6,3\}$ and $\{3,6\}$ are the tilings by squares, regular hexagons and equilateral triangles, respectively.
When $\ineq<$, the graph $\pq$ can be drawn in the hyperbolic plane in a regular way (that is, so that all edges have the same  hyperbolic length and all faces have the same hyperbolic area).
(See Brannan, Esplen and Gray [B, Chapter 6] for hyperbolic length and area.)

Our goal in this paper is to count the vertices of the regular tessellations in the following way.
We choose a vertex as the origin, and partition the vertices into {\it generations\/} according to their distance (as measured by the number of edges in a shortest path) from the origin.
For $n\ge 0$, we denote by $v(n)$ the number of vertices in generation $n$.
Thus, $v(0)=1$, $v(1)=q$, and so forth.
(This sequence does not depend on the choice of the origin, since the automorphism group of a regular tessellation acts transitively on the vertices.
See Coxeter and Moser [C] for a description of the automorphism group.)
The sequence $v(0), v(1), \ldots$ is most easily expressed in terms of its {\it generating function},
$V(z) = \sum_{n\ge 0} v(n)\,z^n$.
For the Platonic solids, the generating function is a polynomial:
$V_{3,3}(z) = 1 + 3z$, $V_{4,3} = 1 + 3z + 3z^2 + z^3$,
$V_{3,4}(z) = 1 + 4z + z^2$,
$V_{5,3}(z) = 1 + 3z + 6z^2 + 6z^3 + 3z^4 + z^5$ and
$V_{3,5}(z) = 1 + 5z + 5z^2 + z^3$.
For the Euclidean tilings, the sequence is (apart from its first term, $v(0)=1$) an arithmetic progression.
It is not hard to see (and will be a consequence of our results) that
$V_{4,4}(z) = 1 + 4z + 8z^2 + 12z^3 + 16z^4 + \cdots = (1+z)^2/(1-z)^2$,
$V_{6,3}(z) = 1 + 3z + 6z^2 + 9z^3 + 12z^4 + \cdots = (1+z+z^2)/(1-z)^2$ and
$V_{3,6}(z) = 1 + 6z + 12z^2 + 18z^3 + 24z^4 + \cdots = (1+4z+z^2)/(1-z)^2$.
For the tilings of the hyperbolic plane, the sequence $v(n)$ grows exponentially with $n$,
but the generating function is still rational in all cases.
In this paper we shall derive the generating function $V_{p,q}(z)$ for all $p\ge 3$ and $q\ge 3$ with $\ineq\le$, thus including those of the Euclidean tilings mentioned above as well as those of the hyperbolic tilings.
For the hyperbolic tessellations $\{4,5\}$, $\{6,4\}$  and $\{3,7\}$  (those tessellations $\pq$ for which $\{p,q-1\}$ is Euclidean), the sequence $v(n)$ can be expressed in terms of the Fibonacci numbers; for other hyperbolic tessellations, the sequence is governed by other, more general, linear recurrences.

The generating functions we seek will be derived from recurrences in a routine way.
(See Wilf [W, Chapters 1 and 2].)
The recurrences, on the other hand, must be derived from combinatorial arguments concerning the tessellations.
We shall divide our work into three sections according to the  combinatorial results needed.
An edge that joins a vertex $a$ in generation $n\ge 0$ to a vertex $b$ in generation $n+1$ will be called a {\it filial\/} edge; we shall say that $a$ is the {\it parent\/} of $b$, and that $b$ is a {\it child\/} of $a$.
In Section 2, we shall deal with the case of $p=\infty$, in which each face is infinite, so there are no finite cycles;  the tessellation is thus an infinite tree in which each vertex has degree $q$.
In this case, every edge is a filial edge; the origin has no parents and $q$ children; and every other vertex has one parent and $q-1$ children.
We easily obtain the result that
$V_{\infty,q}(z) = 1 + qz + q(q-1)z^2 + q(q-1)^2z^3 + q(q-1)^3z^4 + \cdots = qz/(1-(q-1)z)$.

When $p$ is finite, we must treat the cases in which $p$ is even differently from those in which $p$ is odd.
When $p$ is even, every cycle has even length, so the graph is bipartite, and thus two-colorable.
It is clear that assigning one color to the even generations and the other color to the odd generations yields a valid coloring, and thus every edge is a filial edge.
When $p$ is odd, however, it can no longer be true that every edge is a filial edge, for if it were, then coloring the generations according to their parity would yield a valid two-coloring of a graph containing odd cycles, a contradiction.
The edges that join two vertices in the same generation will be called {\it fraternal\/} edges
(if they join vertices, called {\it siblings}, having a common parent) or {\it consortial\/} edges (if they join vertices, called {\it cousins}, having a latest common ancestor in a generation earlier than that of their parents).

In Section 3, we shall deal with the case in which $p\ge 4$ is finite and even.
In this case, it is still true that every edge is a filial edge, but it is now possible for a vertex to have two parents.
A vertex with one parent (and $q-1$ children) will be called a {\it type-A} vertex, while one with two parents (and $q-2$ children) will be called a {\it type-B\/} vertex.
Each face has a {\it latest\/} vertex (which is always a type-B vertex) in some generation $n$, and an {\it earliest\/} vertex (which may be the origin, a type-A vertex or a type-B vertex) in generation
$n-p/2$.

In Section 4, we shall deal with the case in which $p$ is finite and odd.
In this case, it will no longer be true that every edge is a filial edge.
If $p=3$, each vertex other than the origin is joined by {\it fraternal\/} edges to its two 
{\it siblings\/} (the immediately preceding and  following vertices in its generation).
Each such vertex will be either type-A (with one parent, two siblings and $q-3$ children), or
type-B (with two parents, two siblings and $q-4$ children).
There will also now be two kinds of faces, those with a later type-B vertex in generation $n$ and
an earlier fraternal edge joining two siblings (each of which may be either a type-A vertex or a type-B vertex) in generation $n-1$, and those with a later fraternal edge joining two siblings (which each may be either a type-A vertex or a type-B vertex) in generation $n$ and an earlier vertex 
(which may be the origin, a type-A vertex or a type-B vertex) in generation $n-1$.
This case will be reduced to that of an almost-regular tessellation in which every face has degree four.

Finally, if $p\ge 5$ is finite and odd, we have the most complicated situation of all.
In addition to type-A vertices (which again have one parent and $q-1$ children) and type-B vertices (which again have two parents and $q-2$ children), there will now be {\it type-C\/} vertices, each of which is joined by a {\it consortial\/} edge to a single {\it cousin\/} (an immediately preceding or  following vertex in its generation), in addition to having one  parent and $q-2$ children.
There will again be two kinds of faces: those with a latest type-B vertex in generation $n$ and an earliest consortial edge joining two type-C vertices in generation $n-(p-1)/2$, and those with a latest consortial edge joining two type-C vertices in generation $n$ and an earliest vertex
(which may be the origin, a type-A vertex, a type-B vertex or a type-C vertex) in generation
$n-(p-1)/2$.
\sk

\heading{2. Infinite Faces}

For $p=\infty$, the graphs are trees, with $q$ edges meeting each vertex.
The origin has no parents and $q$ children; every other vertex has one parent and $q-1$ children.
Let $a(n)$ denote the number of non-origin vertices in generation $n$.
We have $a(n)=0$ for $n< 1$, and
$$a(n) = (q-1)\,a(n-1) + q\de_{n-1} \eqno(2.1)$$
for $n\ge 1$, where
$$\de_m = \cases{
1, &if $m=0$; \cr
0, &if $m\not=0$. \cr
}$$

Let 
$$A(z) = \sum_{n\ge 0} a(n)\,z^n = qz + \cdots$$
denote the generating function for the number of non-origin vertices in generation $n$.
Multiplying (2.1) by $z^n$ and summing over $n\ge 1$ yields
$$A(z) = (q-1)z\,A(z) + qz,$$
which has the solution
$$\eqalign{
A(z)
&= {qz \over 1 - (q-1)z}. \cr
}$$

If 
$$V(z) = \sum_{n\ge 0} v(n)\,z^n = 1 + qz + \cdots$$
denotes the generating function for the number $v(n)$ of vertices in generation $n$,
then $V(z) = 1 + A(z)$, so we have
$$\eqalign{
V(z)
&= {1+z \over 1 - (q-1)z} \cr
&= 1 + qz + q(q-1)z^2 + \cdots + q(q-1)^{n-1}z^n + \cdots, \cr
}$$
which yields the formula
$$v(n) = \cases{
1, &if $n=0$; \cr
q(q-1)^{n-1}, &if $n\ge 1$. \cr
}$$
\sk

\heading{3. Even Faces}

In this section, we shall treat the case in which $p$, the degree of each face, is finite and even.
Let $p = 2r$, with $r\ge 2$, and suppose $q\ge 3$.
We shall say that a vertex is {\it type-A\/} if it has one parent and $q-1$ children, and that it is of 
{\it type-B\/} if it has two parents and $q-2$ children.
For the tessellation $\{p,q\}$ to be infinite (that is, for $\ineq\le$), we must have $r\ge 3$ or $q\ge 4$, or both.
We claim that under these circumstances, every vertex other than the origin is either type-A or type-B.
(The condition that $r\ge 3$ or $q\ge 4$ is necessary for this claim, since in the cube $\{4,3\}$
the antipode of the origin has three parents and no children, and is thus neither type-A nor type-B.)

Let us first consider the case $q\ge 4$.
Since every vertex has at least one parent, a counterexample would have to have three or more parents.
Suppose, to obtain a contradiction, that such a counterexample $a$ occurs in generation $n$,
and that generation $n$ is the earliest in which such  counterexample can be found.
Let $b$ be a ``middle'' parent of $a$ (that is, a parent $a$ that is neither the leftmost nor the rightmost parent of $a$).
Since $b$ lies in generation $n-1$, it has at most two parents, and thus has at least $q-2\ge 2$ children.
In particular, $b$ has a child $c$ other than $a$.
The child $c$ must lie either to the left or to the right of $a$.
But if $c$ lies to the left of $a$, the edge from $b$ to $c$ crosses the edge from the leftmost parent of $a$ to $a$, while if $c$ lies to the right of $a$, the edge from $b$ to $c$ crosses the edge from the rightmost parent of $a$ to $a$.
Thus we obtain a contradiction to the planarity of the tessellation, proving the claim for $q\ge 4$.

It remains to consider the case $q=3$ and $r\ge 3$.
In this case a counterexample to the claim must take the form of a vertex with three parents and no children.
We shall call such a counterexample a {\it type-I\/} vertex.
A type-B vertex with a type-B parent will be called a {\it type-II\/} vertex.
We shall prove that not only are there no type-I vertices, but there are no type-II vertices.
Suppose then, to obtain a contradiction, that vertex $a$ in generation $n$ is either a type-I or type-II vertex, and that generation $n$ is the earliest generation in which such a counterexample can be found.

Consider first the case in which $a$ is a type-I vertex.
Let $b$, $c$ and $d$ be the parents, from left to right,  of $a$.
Then $b$, $c$ and $d$ are consecutive vertices, for an intervening vertex could not have any children (since the filial edges to such children would cross the filial edges to $a$), and thus would be an earlier type-I vertex in generation $n-1$.
The vertex $c$ cannot have any child other than $a$, since the filial edge to such a child would cross one of the filial edges from $b$ or $d$ to $a$.
Thus $c$ is a type-B vertex.
Let $e$ and $f$ be the left and right parents, respectively of $c$.
Vertex $e$ must be a type-A vertex, else $b$, which lies in generation $n-1$, would be a type-II vertex earlier than the type-I vertex $a$.
Thus $e$ must have exactly one child in addition to $c$.
Since the children of a vertex are consecutive, this additional child must be either $b$ or $d$.
But it cannot be $b$ (since then $a$, $b$, $c$ and $e$ would form a face of degree four, contradicting the hypothesis that $r\ge 3$), and it cannot be $d$ (since then the filial edge from $e$ to $d$ would cross that from $f$ to $c$).
This contradiction shows that $a$ cannot be a type-I vertex.

Consider now the case in which $a$ is a type-II vertex.
Let $b$ and $c$ be the left and right parents, respectively, of $a$.
At least one of these parents is a type-B vertex.
Assume, without loss of generality, that $b$ is a type-B vertex.
Let $d$ and $e$ be the left and right parents, respectively, of $b$.
Vertex $e$ must be a type-A vertex, else $b$, which lies in generation $n-1$, would be a type-II vertex earlier than $a$.
Thus $e$ must have exactly one child in addition to $b$.
This additional child cannot lie to the left of $b$ (since then the filial edge to the child would cross that from $d$ to $b$), and it cannot be the vertex $c$ immediately to the right of $b$ (since then the vertices $a$, $b$, $c$ and $e$ would form a face of degree four, contradicting the assumption that $r\ge 3$).
This contradiction shows that $a$ cannot be a type-II vertex, and thus completes the proof of the claim.

We now proceed to the enumeration of vertices in infinite tessellations $\{p,q\}$ with $p$ finite and even.
For $n\ge 0$, let $a(n)$ and $b(n)$ denote the number of type-A and type-B vertices, respectively, in generation $n$.
We have $a(n)=0$ if $n<1$, and $b(n)=0$ if $n<r$.
For $n\ge 1$, we shall count, in two ways, the number of filial edges between parents in generation $n-1$ and their children in generation $n$.
The origin in generation $0$ has $q$ children in generation $1$.
For $n\ge 2$, each type-A (respectively, type-B) vertex in generation $n-1$ has $q-1$ (respectively, $q-2$) children in generation $n$.
For $n\ge 1$, each type-A (respectively, type-B) vertex in generation $n$ has one (respectively, two) parents in generation $n-1$.
Equating these counts yields
$$a(n) + 2b(n) = (q-1)\,a(n-1) + (q-2)\,b(n-1) + q\de_{n-1} \eqno(3.1)$$
for $n\ge 0$.
For $n\ge r$, each type-B vertex in generation $n$ is the latest vertex of a face whose earliest vertex lies in generation $n-r$.
The origin is the earliest vertex of $q$ faces.
For $n-r\ge 1$, each type-A (respectively, type-B) vertex in generation $n-r$ is the earliest vertex of $q-2$ (respectively, $q-3$) faces (since each face corresponds to a pair of consecutive children of its earliest vertex, and a vertex with $s$ children has $s-1$ paris of consecutive children).
Thus
$$b(n) = (q-2)\,a(n-r) + (q-3)\,b(n-r) + q\de_{n-r} \eqno(3.2)$$
for $n\ge 0$.

Let
$$A(z) = \sum_{n\ge 0} a(n)\,z^n = qz + \cdots$$
and
$$B(z) = \sum_{n\ge 0} b(n)\,z^n = qz^r + \cdots$$
denote the generating functions for the sequences $a(n)$ and $b(n)$, respectively.
Multiplying (3.1) and (3.2) by $z^n$ and summing over $n\ge 0$ yields
$$A(z) + 2B(z) = (q-1)zA(z) + (q-2)zB(z) + qz$$
and
$$B(z) = (q-2)z^rA(z) + (q-3)z^rB(z) + qz^r.$$
Solving these simultaneous equations yields
$$A(z) = {qz(1-2z^{r-1} + z^r) \over 1 - (q-1)z + (q-1)z^r - z^{r+1}}$$
and
$$B(z) = {qz^r(1-z) \over 1 - (q-1)z + (q-1)z^r - z^{r+1}}.$$
If 
$$V(z) = \sum_{n\ge 0} v(n)\,z^n = 1 + qz + \cdots$$
denotes the generating function for the number $v(n)$ of vertices in generation $n$,
then $V(z) = 1 + A(z) + B(z)$, so we have
$$V(z) = {(1+z)(1-z^r) \over 1 - (q-1)z + (q-1)z^r - z^{r+1}}.$$
In this expression, the numerator and the denominator both vanish for $z=1$,
so we may divide both by $1-z$ and obtain
$$V(z) = {(1+z)(1 + z + \cdots + z^{r-2} + z^{r-1}) \over 1 - (q-2)z - \cdots - (q-2)z^{r-1} + z^r}. 
\eqno(3.3)$$
For $\{4,4\}$ ($p=4$, $r=2$ and $q=4$), this expression becomes
$$\eqalign{
V(z)
&= {(1+z)^2 \over (1-z)^2} \cr
&= 1 + 4z + 8z^2 + 12z^3 + 16z^4 + \cdots, \cr
}$$
which yields the formula
$$v(n) = \cases{
1, &if $n=0$; \cr
4n, &if $n\ge 1$. \cr
}$$
For $\{4,5\}$ ($p=4$, $r=2$ and $q=5$), we obtain
$$\eqalign{
V(z) 
&= {(1+z)^2 \over 1 - 3z + z^2} \cr
&= 1 + 5z + 15z^2 + 40z^3 + 105z^4 + \cdots, \cr
}$$
which yields the formula
$$v(n) = \cases{
1, &if $n=0$; \cr
5F_{2n}, &if $n\ge 1$, \cr
}$$
where $F_m$ is the $m$-th Fibonacci number, defined by 
$F_0 = 0$, $F_1 = 1$ and $F_m = F_{m-1} + F_{m-2}$ for $m\ge 2$.


The numerator in (3.3) vanishes for $z=-1$.
If $r$ is odd,  then the denominator also vanishes for $z=-1$, so we may divide both by $1+z$ and obtain
$$V(z) = { 1 + z + \cdots + z^{r-2} + z^{r-1} \over 1 - (q-1)z + z^2 -  \cdots + z^{r-3} - (q-1)z^{r-2} + z^{r-1}}.$$ 
(In the denominator of this expression, the even powers of $z$ have coefficient $1$, while the odd powers have coefficient $-(q-1)$.)
For $\{6,3\}$ ($p=6$, $r=3$ and $q=3$), this expression becomes
$$\eqalign{
V(z)
&= {1+z+z^2 \over (1-z)^2} \cr
&= 1 + 3z + 6z^2 + 9z^3 + 12z^4 + \cdots, \cr
}$$
which yields the formula
$$v(n) = \cases{
1, &if $n=0$; \cr
3n, &if $n\ge 1$. \cr
}$$
For $\{6,4\}$ ($p=6$, $r=3$ and $q=4$), we obtain
$$\eqalign{
V(z) 
&= {1+z+z^2 \over 1 - 3z + z^2} \cr
&= 1 + 4z + 12z^2 + 32z^3 + 84z^4 + \cdots, \cr
}$$
which yields the formula
$$v(n) = \cases{
1, &if $n=0$; \cr
4F_{2n}, &if $n\ge 1$. \cr
}$$
More generally, the denominator of $V(z)$ for $\{6,q\}$ is the same as that for $\{4,q+1\}$.
\sk

\heading{4. Odd faces}

For odd $p\ge 3$, we shall consider two cases in turn:
the case $p=3$ will be reduced to that of $p=4$ treated in the previous section;
the case of odd $p\ge 5$ will be more complicated.

First let $p=3$, and suppose that $q\ge 3$.
For the tessellation $\pq$ to be infinite, we must have $q\ge 6$,
since the tetrahedron $\{3,3\}$, the octahedron $\{3,4\}$ and the icosahedron $\{3,5\}$ are finite.
Every vertex other than the origin is joined by {\it fraternal\/} edges to its two {\it siblings},
the immediately preceding and following vertices in its generation.
These edges link the vertices of each generation into a cycle, with all the cycles surrounding the origin, and the cycles of earlier generations nested within those of later generations.
If these fraternal edges are deleted, the degree of each vertex other than the origin is reduced by two, and each pair of triangular faces that were incident with a given fraternal edge merges into a single quadrilateral face.
In this graph, which we shall denote $\{3,q\}^*$, every face has degree four, and every vertex has degree $q-2$, except for the origin, which has degree $q$.

For $\{3,q\}^*$, the arguments from the preceding section show that every vertex other than the origin is either type-A (which now has one parent and $q-3$ children) or type-B (which now has two parents and $q-4$ children).
Furthermore,  the recurrences for $a(n)$ and $b(n)$ are easily obtained from (3.1) and (3.2), by reducing the coefficients in the homogeneous terms on the right-hand sides by two, while leaving those in the inhomogeneous terms unchanged.
Thus we have
$$a(n) + 2b(n) = (q-3)\,a(n-1) + (q-4)\,b(n-1) + q\de_{n-1}$$
and
$$b(n) = (q-4)\,a(n-2) + (q-5)\,b(n-2) + q\de_{n-2}$$
for $n\ge 0$.
These equations yield
$$A(z) + 2B(z) = (q-3)zA(z) + (q-4)zB(z) + qz$$
and
$$B(z) = (q-4)z^2A(z) + (q-5)z^2B(z) + qz^2$$
for the generating functions $A(z)$ and $B(z)$.
Solving these equations, we obtain
$$A(z) = {qz(1-z) \over 1 - (q-4)z - z^2}$$
and
$$B(z) = {qz^2 \over 1 - (q-4)z - z^2},$$
which yields
$$V(z) = {1 + 4z + z^2 \over 1 - (q-4)z - z^2}$$
for $V(z) = 1+ A(z) + B(z)$.
For $\{3,6\}$ ($p=3$ and $q=6$), this expression becomes
$$\eqalign{
V(z)
&= {1 + 4z + z^2 \over 1 - 2z - z^2} \cr
&= 1 + 6z + 12z^2 + 18z^3 + 24z^4 + \cdots, \cr
}$$
which yields the formula
$$v(n) = \cases{
1, &if $n=0$; \cr
6n, &if $n\ge 1$. \cr
}$$
For $\{3,7\}$ ($p=3$ and $q=7$), we obtain
$$\eqalign{
V(z) 
&= {1 + 4z + z^2 \over 1 - 3z - z^2} \cr
&= 1 + 7z + 21z^2 + 56z^3 + 147z^4 + \cdots, \cr
}$$
which yields the formula
$$v(n) = \cases{
1, &if $n=0$; \cr
7F_{2n}, &if $n\ge 1$. \cr
}$$
We observe that for $p = 3$, $4$ or $5$, the Euclidean cases $\{3,6\}$, $\{4,4\}$ and $\{6,3\}$
yield arithmetic progressions for $v(n)$, while the hyperbolic cases $\{3,7\}$, $\{4,5\}$ and 
$\{6,4\}$ (in which $q$ is just one larger than in the corresponding Euclidean case) yield formulas involving the Fibonacci numbers.

Now let $p=2r+1$, with $r\ge 2$, and suppose that $q\ge 3$.
We shall say that a vertex is {\it type-C\/} if it has one parent, one cousin and $q-2$ children.
For the tessellation $\pq$ to be infinite (that is, for $\ineq\le$), we must have $r\ge 3$ or $q\ge 4$, or both.
We claim that under these circumstances, every vertex other than the origin is either type-A, type-B or type-C.
(The condition that $r\ge 3$ or $q\ge 4$ is necessary for this claim, since in the dodecahedron
$\{5,3\}$, the antipode of the origin  has three parents, and is thus neither type-A, type-B nor type-C.)

This claim is proved by the same arguments as those used to prove the corresponding claim in Section 3.
In particular, when $q\ge 4$, every vertex again has at least two children, and the proof proceeds exactly as in Section 3.
When $q=3$, however, there are many more cases to consider.
A type-I vertex (that is, a counterexample to the claim) could assume any of three forms:
a vertex with three parents, a vertex with two parents and one cousin, or a vertex with one parent and two cousins.
It is again necessary to strengthen the inductive hypothesis to exclude type-II vertices, now defined as type-B or type-C vertices that have a type-B or type-C parent.
A type-II vertex can thus assume any of four forms.
But in every case, a contradiction is reached, either by a violation of planarity or by the existence of a face of degree at most six.
We shall thus omit this tedious consideration of cases, and proceed to use the claim to derive the generating functions for these tessellations.

For $n\ge 0$, let $a(n)$, $b(n)$ and $c(n)$ denote the number of type-A, type-B and type-C vertices, respectively, in generation $n$.
We have $a(n)=0$ if $n<1$, $b(n)=0$ if $n<2r$ and $c(n)=0$ if $n<r$.
For $n\ge 1$, we shall count, in two ways, the number of filial edges between parents in generation $n-1$ and their children in generation $n$.
The origin in generation $0$ has $q$ children in generation $1$.
For $n\ge 2$, each type-A (respectively, type-B, type-C) vertex in generation $n-1$ has $q-1$ (respectively, $q-2$, $q-2$) children in generation $n$.
For $n\ge 1$, each type-A (respectively, type-B, type-C) vertex in generation $n$ has one (respectively, two, one) parents in generation $n-1$.
Equating these counts yields
$$a(n) + 2b(n) + c(n) = (q-1)\,a(n-1) + (q-2)\,b(n-1) + (q-2)\,c(n-1) + q\de_{n-1} \eqno(4.1)$$
for $n\ge 0$.
For $n\ge 2r$, each type-B vertex in generation $n$ is the latest vertex of a face whose earliest edge is a consortial edge that lies in generation $n-r$.
This edge in turn is the latest edge of a face whose earliest vertex lies in generation $n-2r$.
The origin is the earliest vertex of $q$ faces.
For $n-2r\ge 1$, each type-A (respectively, type-B, type-C) vertex in generation $n-2r$ is the earliest vertex of $q-2$ (respectively, $q-3$, $q-3$) faces whose latest edge (which joins two 
type-C vertices) lies in generation $n-r$, and these edges are in turn the earliest edges of $q-2$ (respectively, $q-3$, $q-3$) faces whose latest vertices are type-B vertices in generation $n$.
Thus
$$b(n) = (q-2)\,a(n-2r) + (q-3)\,b(n-2r) + (q-3)\,c(n-2r) + q\de_{n-2r} \eqno(4.2)$$
and
$$c(n) = 2(q-2)\,a(n-r) + 2(q-3)\,b(n-r) + 2(q-3)\,c(n-r) + 2q\de_{n-r} \eqno(4.3)$$
for $n\ge 0$.

Let
$$A(z) = \sum_{n\ge 0} a(n)\,z^n = qz + \cdots,$$
$$B(z) = \sum_{n\ge 0} b(n)\,z^n = qz^{2r} + \cdots$$
and 
$$C(z) = \sum_{n\ge 0} c(n)\,z^n = 2qz^r + \cdots$$
denote the generating functions for the sequences $a(n)$, $b(n)$ and $c(n)$, respectively.
Multiplying (4.1), (4.2) and  (4.3) by $z^n$ and summing over $n\ge 0$ yields
$$A(z) + 2B(z)  + C(z) = (q-1)zA(z) + (q-2)zB(z) +(q-2)zC(z) +  qz,$$
$$B(z) = (q-2)z^{2r}A(z) + (q-3)z^{2r}B(z) + (q-3)z^{2r} + qz^{2r}$$
and
$$C(z) = 2(q-2)z^rA(z) + 2(q-3)z^rB(z) + 2(q-3)z^rC(z) + 2qz^r.$$
Solving these simultaneous equations yields
$$A(z) = {qz(1+z^r)(1-2z^{r-1}+z^r) \over 1 - (q-1)z + 2z^r - 2z^{r+1} + (q-1)z^{2r} - z^{2r+1}},$$
$$B(z) = {qz^{2r}(1-z) \over 1 - (q-1)z + 2z^r - 2z^{r+1} + (q-1)z^{2r} - z^{2r+1}}$$
and
$$C(z) = {2qz^r(1-z) \over 1 - (q-1)z + 2z^r - 2z^{r+1} + (q-1)z^{2r} - z^{2r+1}}.$$
If 
$$V(z) = \sum_{n\ge 0} v(n)\,z^n = 1 + qz + \cdots$$
denotes the generating function for the number $v(n)$ of vertices in generation $n$,
then $V(z) = 1 + A(z) + B(z) + C(z)$, so we have
$$V(z) = {1 + z + 2z^r  - 2z^{r+1} - z^{2r} - z^{2r+1} \over 
1 - (q-1)z + 2z^r - 2z^{r+1} + (q-1)z^{2r} - z^{2r+1}}.$$
In this expression, the numerator and the denominator both vanish for $z=1$,
so we may divide both by $1-z$ and obtain
$$V(z) = {1 + 2z + \cdots + 2z^{r-1} + 4z^r + 2z^{r+1} + \cdots + 2z^{2r-1} + z^{2r} \over 
1 - (q-2)z - \cdots - (q-2)z^{r-1} - (q-4)z^r - (q-2)z^{r+1} - \cdots - (q-2)z^{2r-1} + z^{2r}}. $$
\sk

\heading{5. Conclusion}

We have obtained the generating functions for the number of vertices, classified by generation,
in all of the regular tessellations of the Euclidean or hyperbolic plane.
As a by-product, we have obtained the generating functions for the various ``types'' of vertices in each generation.
It is clear that our methods could be adapted to enumerate edges or faces, or various types of these objects, as well.

We also recall that standard methods can be applied to obtain the asymptotic behavior of the coefficients of a generating function from its analytic behavior (see Wilf [W, Chapter 5], for example).
For $\ineq<$ (that is, when the tessellation is hyperbolic), if 
$V(z) = P(z)/Q(z)$ is the generating function, where the polynomials $P(z)$ and $Q(z)$ each have constant term $1$, then $v(n)\sim Az_0^{-n}$, where $z_0$ is the unique smallest root of $Q(z)$,
and $A = -P(z_0)/z_0 Q'(z_0)$ (the prime denotes differentiation).
We observe that for all the hyperbolic tessellations $Q(1/z) = Q(z)/z^d$, where $d$ is the degree of $Q(z)$, so the smallest root $z_0$ is the reciprocal of the largest root.
\sk

\heading{6. Acknowledgment}

The research reported here was supported
by Grant CCF 0646682 from the National Science Foundation.
\sk

\heading{7. References}

\refbook B; D. A. Brannan, M. F. Esplen and J. J. Gray;
Geometry;
Cambridge University Press, Cambridge, UK,1999.

%
\refbook C; H. S. M. Coxeter and W. O. J. Moser;
Generators and Relations for Discrete Groups;
Springer-Verlag, Berlin, 1980.

\refbook W; H. S. Wilf;
Generatingfunctionology {\rm (third edition)};
A K Peters, Wellesley, MA, 2006.

\bye